\documentclass [12pt, book, reqno, a4paper, czech, chicago, english]{amsart}
\usepackage[dvips]{epsfig}

\usepackage{times}

\usepackage{color}

\newcommand{\R}{I\!\!R}

\def \AA{{\mathcal A}}

\def \UU{{\mathcal U}}

\def \BB{{\mathcal B}}
\def \VV{{\mathcal V}}

\def \R{I\!\!R}

\newcommand{\cqfd} {\hbox {\unskip \kern 6pt \penalty 500\raise -2pt \hbox
{\vrule \vbox to 5pt {\hrule width 4pt \vfill \hrule }\vrule
}\par }}

\newlength\jataille
\newcommand{\figgauche}[3]%
{
\jataille=\textwidth\advance\jataille by -#1
\advance\jataille by -.5cm
\begin{minipage}[a]{#1}
\includegraphics[width=#1]{#2}
\end{minipage}
\vskip2mm
\begin{minipage}[a]{\jataille}
\footnotesize #3 \normalsize
\end{minipage}
}

\usepackage{amsfonts}
\font\tencyr=wncyr10          
\def\cyr{\tencyr
}

\def \ikr{\accent"24 i}

\def\ppref
#1
{$\bullet $
\ref{#1}
}
\newcommand{\lettre}[1]{
\refstepcounter{section}
\vskip 5mm
\centerline{Lettre \thesection}
\sectionmark{Lettre \thesection}
\vskip 3mm }

\usepackage{ulem}
\usepackage[czech]{babel}
\usepackage[cp1250]{inputenc}

\makeindex

\begin{document}

\Large
\centerline{The Axiomatic melting pot}
\centerline{ \it \small Teaching probability theory in Prague during the 1930's}
\vskip 3mm
\centerline{\small \today}
\large
\centerline{
\v{S}t\v{e}p\'anka BILOV\'A\footnote{Pr\'avnick\'a fakulta Masarykovy univerzity, {Brno,  Czech Republic}. bilova@law.muni.cz
 }
, Laurent MAZLIAK\footnote{Laboratoire de Probabilit\'es et Mod\`eles al\'eatoires \& Institut de Math\'ematiques (Histoire des Sciences Math\'ematiques), {Universit\'e Paris  VI,
France}. mazliak@ccr.jussieu.fr}
 and  Pavel \v{S}I\v{S}MA\footnote{Katedra matematiky P\v{r}\'\i rodov\v{e}deck\'e fakulty Masarykovy univerzity, { Brno,  Czech Republic}. sisma@math.muni.cz}
}

\normalsize
\renewcommand{\rightmark}{\uppercase{The Axiomatic melting pot}}
\vskip 5mm
\small
\centerline {\bf  R\'esum\'e}

Dans cet article, nous nous int\'eressons \`a la fa\c con dont la th\'eorie des probabilit\'es a \'et\'e enseign\'ee \`a Prague et en Tch\'ecoslovaquie,  notamment durant les ann\'ees 1930. Nous portons une attention toute particuli\`ere \`a un livre de cours de probabilit\'es, publi\'e \`a Prague par Karel Rychl\'\i{}k en 1938, et qui utilise comme support l'axiomatisation de Kolmogorov, un fait tr\`es exceptionnel avant la deuxi\`eme guerre mondiale.

\vskip 5mm
\centerline {\bf  Abstract}
 In this paper, we are interested in the teaching of probability theory in Prague and Czechoslovakia, in particular 
during the 1930's. We focus specially on a textbook, published in Prague by Karel Rychl\'\i{}k in 1938, which uses Kolmogorov's 
axiomatization, a very exceptional fact before World War II.

\vskip 2mm

\normalsize



\hskip 2cm {\it Secondary} :  28-03, 03E75

\section*{Introduction}

It is now obvious that the first half of the 20th century was a golden age for probability theory. Between the beginning 
of Borel's (1871--1956)  interest in the field in 1905 and the rapid development of the general theory of stochastic 
processes and stochastic calculus in the 1950's, mathematical probability evolved from a rather small topic, with 
interesting but quite scattered results,  to a powerful theory with precise foundations and an active field of 
investigation. There is a huge literature dealing with the description of this burst of interest. Let us just 
mention Von Plato's book \cite{Vonplato1994} whose subject is precisely to describe this new emergence of probability 
and renewal of the theory after  the creations of the 17-18th centuries and the Laplacian impulse of the beginning of the 19th. 
The book contains a huge bibliography in which the reader may find a lot of interesting information to complete the 
picture. The years which follow immediately the First World War are particularly crucial, as they appear to have been 
the moment when the consciousness of the importance and power of set measure theory, intuited by Borel  twenty years ago,  
appeared more clearly to a series of mathematicians, in particular those of the next generation. Among them, the young 
Kolmogorov (1903--1987) made several decisive steps and published the celebrated treatise \cite{Kolmogorov1933} where he 
fixed an axiomatization of probability in the framework of measure theory, an axiomatization which is still used today as 
the basis of the mathematical theory of probability.  It is  too often considered though that this axiomatization was 
accepted immediately and without discussions. On the contrary, the period between 1920 and 1940 was a moment of intense 
debates about the meaning of probability and the possible solutions for giving a reasonable axiomatization to it. The 
passionate discussions in the Berlin and Vienna Circles around Reichenbach and Carnap belong precisely to this moment. 
The reader may consult the excellent paper \cite{ShaferVovk2005} to have a better idea of these debates and a survey of 
the state of probability theory during these years. 

The years following the First World War were also the moment when new born countries, having acquired independence after 
the implosion of the Central Empires and Czarist Russia, were looking for new strategies to improve their scientific 
contacts and their presence in the scientific movement. The example of Czechoslovakia is very interesting in this 
prospect. Born in 1918, the new state paid a particular attention to its cultural and scientific national and 
international life, created two universities in Brno and Bratislava, and intended to intensify its contacts 
with countries such as France and Great-Britain, and at the same time to federate the scientific cooperation with 
the East-European countries. As far as mathematical probability theory is concerned, the role of Bohuslav 
Hostinsk\'y (1884--1951) was determinant. It was studied in particular in \cite{HavlovaMazliakSisma2005} where his 
mathematical and social activity is considered through the prism of his  correspondence with the French mathematician 
Fr\'echet (1878--1973). 

Therefore, Czechoslovakia during the 1920's and 1930's appeared to be a double barycenter, both in the geographical sense 
of the world and in the space of mathematical ideas. Prague was the middle point between Berlin and Vienna and therefore 
a natural point to discuss the diverse interpretations of probability which were set at this moment.  An emblematic 
international conference of philosophy of probability was held in 1929 in Prague where all the tenors were present.
But the country was also the middle point between Paris and Moscow at the precise moment when the nervous center of 
investigation in the field passed from the first to the second. It has already been observed, in particular in 
\cite{Bru2003}, that this particular situation at this particular moment played a part in the interest and decisive 
role of Hostinsk\'y in Markov Chains investigation. 

In this paper, we are interested in another point of view, the question of teaching probability theory. The general 
subject of the appearance of probability theory in the Czech lands has already been considered in several studies such 
as \cite{Macak2005} and \cite{Hyksova2006}. In this paper, we focus on the question of teaching probability in the 
1930's in Prague. 
As already said, the question of axiomatization was very urgent during these years, and Prague was a natural place for 
new educational experiments. We consider this hypothesis through the amazing example of a textbook published by a 
professor of the Czech Technical University in Prague, Karel Rychl\'\i{}k in 1938,  using Kolmogorov axiomatic presentation. 
To the best of our knowledge Rychl\'\i{}k's publication was the only textbook of probability,  with Cramer's one, to use 
Kolmogorov's axiomatization before 1939.  These attempts remained however without posterity and did not survive to the 
war. 

The paper is organized as follows. In a first section, we rapidly expose the state of probability theory teaching in 
the Czech lands until 1938, by presenting its appearance in the syllabus and the main characters who took care of the 
lectures, such as the great mathematician Emanuel Czuber. Then, in a second part, we focus on Czech textbooks on  
probability which appeared in the given period. The third part deals with the biography of Karel Rychl\'\i{}k, 
the fourth is devoted to a rather detailed description of Rychl\'\i{}k's 
textbook, and the fifth one to the way it was received and commented on, in particular in several studies by his follower 
and assistant Otomar Pankraz.

\section{Teaching Probability theory in the Czech lands}

The history of probability in the Czech land has been treated 
in the book {\it V\'yvoj teorie} {\it 
pravd\v e- podobnosti v \v cesk\'ych zem\'\i{}ch do roku 1938
(The Development of Probability Theory in the Czech Lands until 1938)} 
\cite{Macak2005} by Karel Ma\v c\'ak so far. 
This book, published in the series {\it D\v ejiny matematiky (The History of Mathematics)}, 
is not aimed at describing all stages of the development of this mathematical field in great detail but  is devoted to some selected problems. 
It collects Ma\v c\'ak's lectures on the development of 
probability theory which he had read in the previous ten years. 
Quite surprisingly, the book finishes with the very publication of  
Rychl\'\i{}k's textbook \cite{Rychlik1938}, and unlike the present paper does not deal with the following years.  

Probability theory, or more precisely probability calculus, until the  
1920's in the Czech and Slovak lands, was related especially to teaching mathematics at universities, and 
apart from some exceptions there was no original scientific work in this field. 
As the number of mathematics professors at the Prague University was rather 
small,\footnote{In 1882 the Prague University was divided into the Czech and German 
Universities.} and their main role was to give basic lectures on algebra, 
mathematical analysis and geometry, probability theory is to be found in the syllabuses 
of the lectures only exceptionally.\footnote{The first professor to read lectures 
on probability from time to time was 
Wilhelm Matzka (1798--1891) in the 1850's and 1860's.} 
The situation considerably changed by increasing the number of teachers and especially 
by introducing the studies of actuarial mathematics and mathematical statistics after WW I. 
Shortly before this, in 1911, V\'aclav L\'aska (1862--1943) had come to Prague from Lvov and 
had been appointed as Applied Mathematics professor. In 1912 he had already read a one--hour 
public lecture  
 {\it The Introduction to Probability Theory}.\footnote{L\'aska is an author of 
more than 300 works on many areas of applied mathematics, including the textbook 
{\it Probability Calculus} \cite{Laska1921a} and the work {\it Selected Chapters of 
Mathematical Statistics} \cite{Laska1921b} which he certainly used in the above mentioned course.}

A two--year study programme of actuarial mathematics and mathematical statistics 
was organized at the Prague University from the academic year  
1922/23. The leading personality who took care of teaching in both these subject was Emil Schoenbaum 
(1882--1967).\footnote{Schoenbaum studied actuarial mathematics at the University of G\"ottingen in 1906 and following the request of T. G. Masaryk 
he started working in the field of social insurance.  
In 1919 he habilitated for actuarial mathematics and mathematical statistics at the Prague University where he was 
appointed full professor of actuarial mathematics in 1923. 
In 1930 he became one of the initiators of founding the journal 
{\it Aktu\'arsk\'e v\v edy (Actuarial Sciences)}. 
After 1939 he worked abroad where he took part in reforms of social insurance, e.g. in some Latin American countries, 
but also in the U.S.A. and Canada. He died in Mexico in 1967.} 
The lectures on probability theory of this course were read by Milo\v s K\"ossler (1884--1961).\footnote{K\"ossler was 
appointed extraordinary professor of mathematics in 1922. His scientific work focused mainly on number theory 
and he did not publish any results of his own on probability theory.} 

The situation in teaching probability at Austrian technical universities was similar. 
Until the end of the 19th century there were, normally, no special lectures given on this subject 
 and they were read only as optional lectures by private 
docents.\footnote{In the middle of the 19th century there existed special study programmes 
at technical universities for a certain period where ''political arithmetic'' was 
taught in which probability and statistics played a significant role.} 
One of them was Emanuel Czuber (1851--1925)\footnote{Emanuel Czuber was born on 19th January, 1851 in Prague where he graduated at the 
German Technical University in 1874. Czuber, with original surname \v Cubr, 
came from a completely Czech family, nevertheless he obtained his education at German 
schools. When he started teaching in Vienna, he used only the German language for all 
conversations. In 1872--1875 he was an assistant of geodesy professor  
Ko\v ristka. Then he became a mathematics and descriptive geometry teacher 
at a Prague grammar school. In 1876 he habilitated and as a private docent he 
read two--hour lectures on probability theory, 
the least square method and mathematical statistics. 
In 1886--1891 Czuber was mathematics professor at 
the Technical University in Brno, then at the Technical University in Vienna until 
his retirement in 1921. He died on 22nd August 1925 in Gnigl at 
Salzburg. More detailed information on the life and scientific work of 
Emanuel Czuber, including the list of his publication, can be found e.g. in 
\cite{Dolezal1928}.
Nevertheless a detailed analysis of Czuber's scientific life, his position 
within Austrian and German mathematical community has not been done yet.
During his life Czuber was given a number of significant functions in the area of 
education as well as actuary.  It was him who initiated actuarial technical 
studies in Austro--Hungary, -- details will be given in the following 
part. In 1898--1900 he was the president of the Association of 
Austrian--Hungarian actuarial technicians; from 1900 the chair of the Association of 
Austrian--Hungarian insurance companies; in 1909 he was the chair of the VIth International 
Congress for Actuarial Sciences in Vienna; 
in 1890/91 he was the Rector of the Technical University in Brno, in 1894/95 the Rector 
of the Technical University in Vienna; in 1876--1886 the editor of {\it
Technische Bl\"atter}; in 1897--1921 the editor in chief of  {\it Zeitschrift f\"ur 
Realschulwesen}; the Chair of the Austrian section of the International 
Committee for teaching mathematics, etc. 
His scientific work includes especially 
works on probability theory, 
the theory of errors and adjustment,
geodesy and actuary. 
Shortly after his habilitation he translated 
lectures on probability theory of A. Meyer from Lutych \cite{Czuber1879} into German 
and five years later he wrote the first probability textbook of his own {\it
Geometrische Wahrscheinlichkeiten und Mittelwerte} \cite{Czuber1884} which was published 
also in French. During his stay in Brno he prepared the book \cite{Czuber1891}.}
 who became a private docent in 1876 when he habilitated 
in the field of the theory and practice of 
adjustment. Later he extended his {\it veniae docendi} to probability theory.

In 1903 the first edition of his book {\it Die Wahrscheinlichkeitsrechnung
und ihre Anwendung auf Fehlerausgleichung, Statistik und
Lebensversicherung} \cite{Czuber1903} was published, a huge work which was 
later divided into two volumes (altogether more than 900 pages), the part on probability 
was published in the sixth edition as late as 1941. Thus, this book influenced 
teaching probability not only in German speaking countries for nearly 40 years, 
it was a basic textbook for Czech mathematicians for a long time as well. We may suppose 
that it was used mainly for preparation of future workers in actuary.\footnote{Czuber 
published two other books in 1920's, \cite{Czuber1921,Czuber1923b}.}

Czuber was engaged also in historical and philosophical questions 
of probability theory. In 1899 he published a nearly 300 pages study \cite{Czuber1899}
in {\it Jahresbericht der Deutschen
Mathematiker-Vereinigung} where he described the development of probability theory 
in the 18th and 19th centuries.\footnote{The list of literature in the study includes around 
500 works and it can be considered the second most significant work on 
the history of probability theory, the first being Todhunter's book 
\cite{Todhunter1865}. Czuber also contributed to {\it Encyklop\"adie der 
mathematischen Wissenschaften} with the paper \cite{Czuber1900}.} 
Two years before his death the book {\it Die
philosophischen Grundlagen der Wahrscheinlichkeitsrechnung}
\cite{Czuber1923a} was published.

Private docent Augustin P\'anek (1843--1908) started reading lectures on probability theory at the 
Czech Technical University of Prague in the 1870's.\footnote{In 1872 Augustin P\'anek (1843--1908) habilitated for 
the field of integrals, nevertheless 
he gave mainly lectures on probability theory and the least squares method until 
his appointing an extraordinary professor in 1896. Detailed information about 
Augustin P\'anek can be found in the work of M. Be\v cv\'a\v rov\'a \cite{Becvarova2004}.}  
He read lectures on probability theory until his death, which means even in the period 
when the actuarial technical course was open. The syllabus of P\'anek's 
two--hour lectures 
can be found in the list of lectures at the Czech Technical University: 

\begin{quote}
Absolute, relative and complex probability. Geometric probability. Bernoulli and Poisson Theorems. 
Objective and subjective expectations. Probability a posteriori. Bayes Formula. Laplace Theorem. 
On insurance. Probability and judgment. Historical overview of probability calculus and the least squares 
method.\footnote{Absolutn\'\i{}, relativn\'\i{} a slo\v zit\'a pravd\v epodobnost. Geometrie pravd\v epodobnosti. V\v eta
Bernoulliho a Poissonova. Objektivn\'\i{} a subjektivn\'\i{} nad\v eje. Pravd\v epodobnost
a posteriori. Pravidlo Bayesovo. Theorem Laplace\accent23uv. O poji\v s\v tov\'an\'\i{}.
Pravd\v epodobnost o sezn\'an\'\i{} sv\v edk\accent23u. D\v ejinn\'y n\'astin po\v ctu pravd\v epodobnosti
a methody nejmen\v s\'\i{}ch \v ctverc\accent23u.} 
\end{quote}

P\'anek never attempted to write a textbook on probability theory, but his works include 
articles on probability theory. Those, more or less popularizing, papers 
published in {\it \v Casopis pro
p\v estov\'an\'\i{} matematiky a fysiky} might 
indicate what 
P\'anek was teaching. 
We do not find there, apart from some exceptions, any general considerations, they include  
in fact solving some simple or more difficult probability exercises. 

An important impulse for the development of probability theory teaching came from 
opening the actuarial technical courses at the Austrian--Hungarian technical universities 
at the turn of the century. Those courses existed in Austria as well as at Czech and 
German technical universities in Czechoslovakia until the end of 1930's. 

The first such course was established by Emanuel Czuber at the Vienna Technical University 
in the academic year 1894/95. A very important part of 
the courses was covered by mathematical subjects. Apart from basic mathematics lectures, common to other study programmes, 
the students attended lectures on actuarial mathematics, mathematical statistics, and 
probability theory. The course was open as a three--year programme, but it was 
shortened to two years already in 1897. It became a model for other schools and despite considerable efforts 
to provide it as a four--year study programme equal to other courses at technical university, 
it did not happen so. 

Probability theory at the Vienna Technical University was read by Emanuel
Czuber from 1891 until the end of WW I.\footnote{In 1920 a position of probability 
theory professor was established which was not, however, occupied and until 
the end of the war it was only substituted.} 
His lecture took two hours, after a three--hours opening of the course. It was included into 
the second year and had the following syllabus:

\begin{quote}

Notion of probability. Direct determination of probability.
Indirect determination of probability. Repeted trials - theorems of  Bernoulli and Poisson. Mathematical expectation, 
mathematical risk. Probability of causes and of future events on the basis of experience. Bases of error theory. 
Least squares method. Theory of collective measure.\footnote{Wahrscheinlichkeitsbegriff. Direkte Wahrscheinlichkeitsbestimmung. Indirekte
Wahrscheinlichkeitsbestimmung. Wiederholte Versuche --- die S\"atze von
Bernoulli und Poisson. Mathematische Hoffnung und mathematische Risiko.
Wahrscheinlichkeit von Ursachen und k\"unftigen Ereignissen auf Grund der
Erfahrung. Elemente der Fehlertheorie. Methode der kleinsten Quadrate.
Kollektivma{{\ss}}lehre.}

\end{quote}

Mathematical statistics which was taught three hours during the whole course was 
taught by Ernst Blaschke (1856--1926), the author of the textbook \cite{Blaschke1906},  
in 1896--1926. Actuarial mathematics was read by 
Viktor Sersavy (1848--1901) until 1901 and from 1901 to 
1938 by Alfred Tauber (1866--1942).\footnote{More detailed information about 
actuarial mathematics teaching in Vienna can be found in 
\cite{Rybarz1965,Ottowitz1992}.}

The second school where the course was established, in 1904, was 
the Prague Czech Technical University. It was initiated by Gabriel Bla\v zek
(1842--1910) who had been reading lectures on actuary mathematics 
since 1901. The lectures were then taken over by Josef Bene\v s (1859--1927)
who was appointed a professor in this field after WW I and was also teaching 
mathematical statistics. After Bene\v s's death actuary mathematics teaching was 
taken over by Jaroslav Janko (1893--1965).

Probability was first taught two hours the whole second year, then four hours in the 
winter semester only. Probability theory was not examined at the final examination, students 
only had to prove their knowledge during the studies. 
The two-hour lectures were given by professor
P\'anek until his death, then by his successor Franti\v sek Vel\'\i{}sek (1877--1914) who read lectures 
on probability four hours in winter semester. After Vel\'\i{}sek died on the front in 1914, 
his teaching was taken over by Karel Rychl\'\i{}k whose duties included lecturing on 
probability theory even after he was appointed professor in 1920. 
From the academic year 1933/34 Rychl\'\i{}k gave two--hour lectures 
 {\it Introduction to probability calculus and mathematical statistics}
in summer semester, however, this lecture was intended for students of 
other programmes. Unfortunately, the content of his lectures is not known. 

Mathematical statistics was first lectured as an independent subject three hours during the 
whole second year, then it became part of actuarial mathematics lectures, and afterwards an 
independent subject again. The syllabus at the beginning of the 1930's was the following:

\begin{quote}
The history of statistics. The methods of statistical investigation. Describing units from 
the point of view of qualitative and quantitative indicators. Measures of variance. The theory of index numbers. 
Mortality tables and their construction. Intensive and time measures, their 
calculation from a given material. Foundations of adjustment theory. The methods of statistical research 
on causal relations. The stability of statistic numbers. Stochastic dependence between 
qualitative and quantitative indicators. The analysis of time series.\footnote{D\v ejiny statistiky, Technika
statistick\'eho \v set\v ren\'\i{}, Popisov\'an\'\i{} soubor\accent23u z hlediska kvalitativn\'\i{}ch a
kvantitativn\'\i{}ch znak\accent23u. M\'\i{}ry rozptylu. Teorie indexn\'\i{}ch \v c\'\i{}sel. Tabulky
\'umrtnosti  a jejich konstrukce. M\'\i{}ry intenzivn\'\i{} a \v casov\'e; jejich v\'ypo\v cet z
dan\'eho materi\'alu. Z\'aklady po\v ctu vyrovn\'avac\'\i{}ho. Metody statistick\'eho b\'ad\'an\'\i{}
o p\v r\'\i{}\v cinn\'ych spojen\'\i{}ch. Stabilita statistick\'ych \v c\'\i{}sel. Stochastick\'a
z\'avislost mezi kvalitativn\'\i{}mi znaky a mezi kvantitativn\'\i{}mi znaky. Rozbor
\v casov\'ych \v rad.}
\end{quote}

Very shortly after having introduced actuarial technical courses at the Czech Technical 
University similar course was opened also at the German Technical University in Prague. It was in 
1906 and actuarial mathematics was taught by 
Gustav Rosmanith (1865--?) until the mid 1930's. He taught also mathematical statistics 
which was taken over by his successor Josef
Fuhrich (1897--?). Two--hour lectures on probability theory were read by Karl Carda
(1870--1943). The syllabus of his lectures at the end of the 1920's was the following:

\begin{quote} 

Problems about urns. Classical problems: problem of  de Moivre die, meeting problem, sharing problem. Repeated trials, 
J. Bernoulli's theorem. Poisson approximation formula. 
Probability of causes. Bayes theorem.\footnote{
Urnenaufgaben. Klassische Probleme: W\"urfelproblem v. Moivre. Rencontreprobl.
Teilungsproblem. Wiederholte Versuche, Theorem v. Jakob Bernoulli.
N\"aherungsformel v. Poisson. Wahrscheinlichk. von Ursachen. Theorem v.
Thomas Bayes.}
\end{quote}

During the war, lectures on probability theory were given by Fuhrich who 
taught actuarial mathematics as well as statistics. The lectures were common for 
the Technical University students and the German University students. 

The last technical university in the area of Czechoslovakia in which 
the actuarial technical course was offered, in 1908, was the German Technical University 
in Brno.\footnote{Details about mathematics teaching at the German Technical University 
in Brno and at universities in Prague are dealt with in \cite{Sisma2002}.
Such course was never introduced at the Czech Technical University in Brno.} There had been 
no special lectures on probability theory or mathematical statistics given at that 
university until the end of 19th century. In 1906 Friedrich Benze (1873--1940), an assistant of mathematics, was charged 
with giving such lectures as a 
paid docent.  
Benze, who did not publish any work was teaching probability and 
statistics until 1939. 
The syllabus of his lectures remained the same all that time, in the academic year 1914/15 
it was the following:

\begin{quote}

Probability Theory I, winter semester 2/0:

\end{quote}

\begin{quote}

Mixing, divisions and
composition of finite sets. Mean determination. Quotas, argument, intersection and dispersion (?). Distributions in 
Bernoulli scheme, Poisson scheme, Compounded Poisson scheme.\footnote{
Unordnungen,  Zerlegungen  und  Zusammensetzungen  endlicher
Mengen. Mittelbildung.  Quotenmittel, Argument, Durchschnitt
und  Streuung.  Verteilungsgesetze  im  Bernoullischen,  im
einfachen und zusammengesetzten Poissonschen Schema.}

\end{quote}

\begin{quote}

Probability Theory II (and applications), summer semester 3/0:

\end{quote}

\begin{quote}

The urn model and its numerical elements:  Bernoulli scheme, 
Poisson and generalized Poisson urn models. The urn model. Bayes theorem. Mathematical risk theory. Theory of fair prices. Comparison between the theoretical urn model and the empirical distribution  : theory of errors,
mathematical theory of mass phenomena.\footnote{Das    Urnenschema und seine numerische Elemente: Bernoullische,
einfache und verallgemeinerte Poissonsche Urnenschema. Das verbundene
Urnenschema. Bayessche  Satz. Theorie des mathematischen Risikos. Theorie
der Wertgleichungen. Vergleich  des theoretischen  Urnenschemas mit
beobachteten Verteilungen: Fehlertheorie, Mathematische Theorie der
Massenerscheinungen.}

\end{quote}

\begin{quote}

Mathematical Statistics, winter semester and summer semester 2/0:

\end{quote}

\begin{quote}

Theory of collective measure. Formal theory of statistical series. Formal theory of population. Biometric functions.
Mixing of statistical distributions. Immigration and emigration. Statistical theory 
of dispersion of sets. Building of statistical tables: Interpolation and compensation of mortality and disability tables. 
Local and analytical compensation.\footnote{Kollektivmasslehre. Formale Theorie statistischer Reihen. Formale
Bev\"olkerungstheorie. Biometrische Funktionen. Mischung von statistischen
Verteilungen. Ein- und Auswanderung. Dispersionstheorie statistischer
Mengen. Konstruktion statistischer Tafeln: Interpolation und Ausgleichung
von  Sterblichkeits-  und Invalidit\"atstafeln. Lokale und analytische
Ausgleichung.}

\end{quote}

Actuarial mathematics was taught by Ernst Fanta
(1878--1939), an actuarial mathematician from Vienna, in 1906--1920. He had worked 
in the seminar of Georg Bohlmann (1869--1928) in G\"ottingen in 1901--1902. During the 1920's
actuarial mathematics was taught by Ferdinand Schnitzler (1857--1933), and then, 
until the war by Oskar Kubelka (1889--?). Despite great efforts of the teaching staff, 
the professorship of actuary mathematics was never established in Brno.

\section{Textbooks}

The 1860's was the time of a fast development of both secondary school 
and university education in Czech language in the Czech lands. 
In 1869 the Prague 
Technical University was divided into the German and 
Czech Technical Universities, and a similar division happened within the Prague University 
in 1882. The first Czech mathematics textbooks were written for the secondary schools 
to replace the existing German books. A significant role was played by the Union 
of Czech Mathematicians and Physicists from the very beginning as they paid attention to 
publishing secondary schools as well as university textbooks.The first task was naturally to cover the area of basic mathematical fields which meant the basis 
for the education of secondary school teachers, then engineers. 
As probability theory was not among those subjects, we cannot find any Czech book 
devoted to this topic until 1921. The situation was similar in other fields of higher mathematics 
in which German textbooks continued to be used and the specialists studied also the literature 
in other foreign languages. 

The first Czech textbook on probability theory was the book {\it Po\v cet
pravd\v epodobnosti (Probability Calculus)} \cite{Laska1921a} from 1921. This book of L\'aska's 
has only 128 pages, including appendices, which are divided into three chapters:
a priori Probability , a posteriori Probability , and Geometric probabilities. 
As E. Schoenbaum stressed in quite a critical review \cite{Schoenbaum1922}
in {\it \v Casopis} ``{\it a characteristic feature of the book is the author's interest in 
philosophical and noetic side of the calculus\footnote{The footnotes prove L\'aska's good 
knowledge of older as well as more recent literature which is devoted to the development 
and philosophical questions of probability theory.} and beside this he uses 
symbolic algorithms in an unusual scope; an undisputable advantage of 
this method is, however, balanced by the loss of space necessary for deducing 
the rules for calculations with the symbols for which the author introduces 
new, sometimes perhaps too complicated signs.}'' L\'aska mentions also the need for 
axiomatization of probability theory in the introduction:

\begin{quote}
{\it From the point of view of pure mathematics 
our considerations on probability theory should in fact begin in a similar 
way as Hilbert starts his geometry. ``Let us think about a concept which 
we name {\it mathematical probability}. We do not know what its transient 
meaning is, nor we have to or need to know either. 
Indeed, it would not be good if we wanted to know. All that is necessary 
to  know about the concept will be told by the axioms.''
Unfortunately the axiomatics of probability calculus has not been worked out  
in a way to be able to present the introduction to the theory.}
\end{quote}

L\'aska mentions the axiomatization attempts of Borel, Broggi, 
Bohlmann\footnote{Broggi and Bohlmann will also appear in Rychl\'\i{}k's 
textbook - see Part 4.} 
and  gives a short exposition of von Mises collectives theory of 1919, 
insisting on the second Mises' axiom as a requirement for a ''perfect mixing'' of the sequence.  
{\it The main obstacle of Mises theory lies in the mathematical formulation of the
perfect mixing, and this, as it appears, and as Mises's more than 
extensive works prove, is not an easy matter.} In another place, L\'aska writes that 
{\it Mises' mathematical theory is in fact an analysis of geometric 
probabilities. This is its theoretical value, but also a mistake in 
the sense of application.}

L\'aska's book was published in the edition of the Czech Technical Matrix (Organization)  
in Prague, not by the Union of Czech Mathematicians and Physicists. 
There was an agreement between these two companies that, due to 
a limited number of Czech readers,  they would not 
compete in publishing books on the same topic. 
That was the reason why in 1925 Hostinsk\'y's request for publishing his lectures on 
probability theory given in Masaryk University of Brno at the beginning of 1920's 
was refused by the Union with an explanation that his book would in fact deal 
with the same matters as L\'aska's textbook, with the exception of geometric 
probability. Hostinsk\'y was invited to write a booklet on geometric 
probability, which he did, and in 1926 his book {\it Geometric 
Probabilities} \cite{Hostinsky1926a} was published.\footnote{Hostinsk\'y made an effort 
to publish his lectures in litography, but the Union refused. See 
the Masaryk University Archives, Bohuslav Hostinsk\'y Fond.}

Hostinsk\'y wrote in the introduction:

\begin{quote}
{\it Various objections have been raised against the notion of geometric 
probability which is so important for physics (e.g. in kinetic theory of 
gases). Today, however, they are mostly of historical importance as the last decades 
have critically explained, mainly due to Borel, the concept of probability and 
its relation to physical applications.}

{\it This book has two purposes. First it gives the basic theorems on geometric 
probabilities and deals with exercises which are interesting from 
purely geometric point of view; a special chapter is devoted to considerations on 
attempts which can approximately confirm theoretical formulas for probabilities.}
\end{quote}

Hostinsk\'y points out to a much more extensive book by Czuber
\cite{Czuber1884} which he quotes on several places. He emphasizes that 
he uses for solutions some special exercises of Poincar\'e's ``method
of arbitrary functions'' which is still not sufficiently known among Czechoslovak 
mathematicians.\footnote{Poincar\'e's method of arbitrary functions is introduced 
on pp. 70--85 where the reader is also acquainted with Hostinsk\'y's results of 
his work on a new solution of Buffon's exercise about the needle \cite{Hostinsky1917}
or his work \cite{Hostinsky1926b}. See also \cite{HavlovaMazliakSisma2005}.}

Hostinsk\'y expresses his thanks to his assistant Josef Kauck\'y for the help 
with correcting the text in the introduction to the book. Kauck\'y belonged to 
Bohuslav Hostinsk\'y's pupils who were interested in the questions of 
theoretical physics as well as probability theory.\footnote{Josef Kauck\'y (1895--1982) 
habilitated for mathematics at the University of Brno in 1928, he spent 
several years teaching at a secondary school in Brno, and in 1938 
he was appointed mathematics professor at the Technical University 
 of Ko\v sice in Slovakia which was transferred to Bratislava during the war. 
From 1946 he was teaching at the Technical University in Brno, and later 
at the Military Academy. He is known for his book {\it Kombinatorick\'e
identity (Combinatorial Identities)} within the Czechoslovak mathematical 
community. He was involved in the areas of difference equations, 
projective and differential geometry, and he also wrote several books on 
probability theory.}
In 1934 he published the book {\it Introduction to Probability Calculus and 
Statistical Theory} \cite{Kaucky1934}. The book was published from the initiative 
of electrotechnology professor V. List who belonged to the most important figures 
of Czechoslovak technical education between the wars, with the support of 
the Czechoslovak Electrotechnology Union. 
It is a small booklet (78 pages) in which Kauck\'y included besides the presentation 
of the basic concepts of probability theory\footnote{Also Kauck\'y presented the method 
of arbitrary functions in his book.} and mathematical statistics, also applications
in the theory of 
errors and statistical mechanics. As Kauck\'y writes in the introduction,
{\it The small scope of the book naturally resulted in the fact that I limited myself 
to the most elementary considerations at some places. Thoughts of 
philosophical nature and historical comments were also omitted.}

Another book which is interesting from the point of view of probability theory teaching in the Czech lands 
is {\it Z\'aklady teorie statistick\'e metody (The Basics of Statistical Method Theory)} \cite{Kohn1929}
written by Stanislav Kohn (1888--1933), a private statistics docent at the Russian Faculty of Law
in Prague. Kohn's book is based on his lectures given in 
Tiflis, Paris and Prague and it is dedicated to the memory of Kohn's teacher A.A.
\v Cuprov. The book is nearly 500 pages long and a short analysis of approximately 
100 pages on probability theory can be found in Ma\v c\'ak's book \cite{Macak2005} (pp. 110--113). 
Kohn's book includes a big number of historical and especially philosophical comments, 
which makes it quite unique in the Czech literature. The list of the used and recommended 
literature is also very extensive, it comprises of more than 30 pages. 

Another book, {\it Probability Calculus} \cite{Hostinsky1950},  was written by  Bohuslav Hostinsk\'y (1884--1951) in 1950, based on lectures given in the 1930's. In the introduction Hostinsk\'y points out that the lectures read in the 1930's differed considerably 
from those given in the 1920's. His expositions treat also the parts of his own scientific research. 
The first part includes a chapter on geometric probability and the second part consists mainly of the 
theory of Markov chains. The application in physics is not present there, Hostinsk\'y intended to deal with this 
topic in another book. This one, however, was not finished as Hostinsk\'y died in 1951. 
Hostinsk\'y's book does not include axiomatic building of probability theory, the reader is acquainted 
with the classical definition and then the author builds probability theory on the basis of 
theorems on addition and multiplication of probabilities which he calls axioms. 

The only textbook, where an axiomatic construction is presented,  seems to have been the one written in 1938 by 
Karel Rychl\'\i{}k, a professor at the Czech Technical University in Prague whom we shall present now.

\section{Karel Rychl\'\i{}k: biography}

Karel Rychl\'\i{}k was born on 16 August 1885. His life is 
described in great detail in the book of M. Hyk\v sov\'a \cite{Hyksova2003} which is the author's Ph.D. thesis 
at Charles University in Prague. The information given in this 
part of the paper is mostly based on that book. 

Rychl\'\i{}k started his secondary school education in 1896 by attending 
grammar school in Chrudim, continuing in his native town Bene\v sov from 1897, and finishing 
with a graduation exam at the Academic Grammar School in Prague in 1904 where the family had moved in 1900.\footnote{Two 
years before Bohuslav Hostinsk\'y  had graduated at the same grammar school.} 
Rychl\'\i{}k had already showed a great interest in mathematics as a secondary school student and his name could be often found 
among those solving successfully the exercises for students published in the {\it \v Casopis pro p\v estov\'an\'\i{} matematiky a fysiky}. 
His mathematics teachers at the prestigious Academic Grammar School included Jan Vojt\v ech
(1879--1953), his later colleague at the Czech Technical University in Prague. 

Rychl\'\i{}k commenced his university studies in the winter semester of 1904/05, he studied at the Faculty of Arts of the 
Czech University in Prague which was then called Charles-Ferdinand University. It was the period of improved  
mathematics teaching at the university as the positions which had been vacated by the deaths of two old and ill professors, 
of professors Franti\v sek Josef Studni\v cka (1836--1903) and Eduard Weyr (1852--1903), were  being occupied by
Karel Petr (1868--1950) and Jan Sobotka (1862--1931). It was especially the algebraist Petr who mostly influenced 
Rychl\'\i{}k's scientific work. Rychl\'\i{}k became a member of professor Petr's seminar, and he 
also gave lectures on his seminar work at the meetings of the Union of Czech Mathematicians
and Physicist whose member he had become after starting the university. He received a scholarship 
of Bernard Bolzano Foundation for his outstanding study results in 1906. 

In the academic year 1907/08 Rychl\'\i{}k obtained a national scholarship and studied 
at the Facult\'e des Sciences in Paris. He attended, apart from others, lectures of Jacques
Hadamard, Emile Picard or Gaston Darboux.\footnote{A year later, it was 
Bohuslav Hostinsk\'y who studied in Paris.} At the Coll\`ege de France 
he listened to lectures on number theory which were given by 
Georges Humbert. During his Paris stay Rychl\'\i{}k already worked on his thesis which he submitted at the Prague 
University in 1908. 
A part of the thesis was published in the {\it \v Casopis pro p\v estov\'an\'\i{} matematiky
a fysiky} with title {\it On a Group of  order 360 (O grup\v e \v r\'adu 360)}. 
After a successful {\it rigorosum examination} Rychl\'\i{}k was awarded the degree 
Doctor of Philosophy on 30 March 1909. Before that,  he had passed also the 
{\it teacher examination} for teaching mathematics and physics at secondary schools and 
he did his teacher practice year at grammar school in \v Zitn\'a ulice, Prague.  

From the beginning of 1909 Rychl\'\i{}k was helping in the mathematical seminar of the Faculty 
of Arts as a non-paid assistant. He worked as a paid assistant from October 1909 till the end of 
June 1913 when he became an assistant at the department of mathematics with 
professor F. Vel\'\i{}sek (1877--1914) at the Technical University.\footnote{Karel Rychl\'\i{}k succeeded his brother 
Vil\'em in that position.  Vil\'em Rychl\'\i{}k (1887--1913) was a talented mathematician who unfortunately
died very young, at the age of twenty six.}

In November 1910 Rychl\'\i{}k submitted at the university an application for granting {\it veniae docendi} 
from mathematics with a habilitation work {\it A Contribution to the theory of forms (P\v r\'\i{}sp\v evek
k teorii forem)}. Probably due to his quite small publication 
activity it was only in March 1912 that his habilitation procedure was successfully finished 
and Rychl\'\i{}k was appointed a {\it private docent} (university lecturer) of mathematics. 
In the academic year
1912/13 he started reading lectures at the university. He stopped with regular lecturing only in 
1925 when he was a professor at the Technical University. 
Rychl\'\i{}k continued giving irregular lectures at Charles University until World War II.
Let us remind that in the years 1912--20 another private docent working at the university was 
Bohuslav Hostinsk\'y. During the World War I they alternated (together with Bohumil Byd\v zovsk\'y (1880--1969)) 
in giving introductory lectures on differential and integral calculus and analytical geometry in space. 
As a private docent Rychl\'\i{}k usually read lectures on algebra 
and number theory topics, which were the fields of his own mathematical research. 

When Rychl\'\i{}k had come to the Technical University he 
applied for transferring his habilitiation from the university and he became 
a private docent at the Technical University too. While professor 
Vel\'\i{}sek was immediately after the beginning of World War I called up 
and died in the war, Rychl\'\i{}k was never called up and during the war substituted all 
Vel\'\i{}sek's basic lectures including the theory of probability. 
It seems very likely that had it not been for the given reasons, 
Rychl\'\i{}k would have not found a way to probability theory. 

The formation of Czechoslovakia brought significant changes also in university 
education. The time came to fulfilling a long lasting attempts to create the second Czech 
university when in 1919 Masaryk University in Brno was founded. The Czech Technical University 
in Prague was re--organized. These two events resulted in increasing the number of mathematics professor
positions. Hostinsk\'y, as well as Rychl\'\i{}k were suggested to be appointed  
professors by the board of professors at Charles University in March 1919, however
the appointment was never put into practice. 
Hostinsk\'y was appointed full professor of theoretical physics at Masaryk University 
in 1920, and Rychl\'\i{}k was appointed extraordinary professor at the Prague Technical University 
in the same year. He was appointed full professor there in 1924. 

Rychl\'\i{}k worked at the Technical University until the closure of 
Czech universities by the Nazis in November 1939.\footnote{In the academic year 1934/35 he was 
the Dean of the Faculty of Mechanical and Electrical Engineering.} 
During the whole period he gave lectures especially on differential and 
integral calculus, and those parts of high mathematics necessary for future engineers. 
His scientific interests would certainly suit a university professor post, however, at Charles University 
the number of teachers grew very slowly around 1920, and algebra and number theory, 
the two main fields of Rychl\'\i{}k's interest, were covered by Karel Petr until WW II. 

After the closure of Czech universities, like all other professors Rychl\'\i{}k was 
sent to the so--called expectant leave which meant his existing salary. 
The possibilities of scientific work were limited during the war mainly for 
an extremely restricted approach to literature, new as well as old -- locked in university libraries. 
Rychl\'\i{}k prepared for publishing the second edition of his textbook 
on elementary number theory and in 1944 started translating Glivenko's 
textbook on probability theory \cite{Glivenko1939}.

Rychl\'\i{}k did not return to the Technical University after the war. 
He had been accused and in May 1946 found guilty of 
{\it breaching the loyalty to the Czechoslovak Republic and transgressing against 
the national honour} by the cleansing committee of the Provincial National Committee in Prague. 
Based on this conviction Rychl\'\i{}k was punished by being made permanently retired with 
a pension reduced by 30\%. 
Rychl\'\i{}k appealed to the cleansing committee at the Ministry of Education which, 
however, rejected his appeal at the end of 1947. The reasons why Rychl\'\i{}k was made retired 
are given by Hyk\v sov\'a in \cite{Hyksova2003} -- e.g. according to the opinion of Vladim\'\i{}r
Ko\v r\'\i{}nek (1899--1981), Rychl\'\i{}k's assistant in 1927--1931 and later professor at 
Charles University, {\it approximately from the beginning of the civil war in Spain, 
Karel Rychl\'\i{}k began to show fascistic opinions and his unconcealed admiration for Fascist 
regimes. This even increased during the so--called second republic.}\footnote{This was 
the period from 1 October 1938 until 15 March 1939.
Ko\v r\'\i{}nek also states that he later reduced his contacts with Rychl\'\i{}k, and so he cannot say 
anything about his opinions during the war. It is probably an attempt to 
save Rychl\'\i{}k from a harsher punishment.}
Other witnesses in the trial with Rychl\'\i{}k pointed out that Rychl\'\i{}k
approved of many measures of German organisation and administration, and that he
manifested his belief that the Germans were going to win the war. According 
to the committee, Rychl\'\i{}k must have been aware that such words were unsuitable to be expressed 
by a university professor in public 
in the time of persecution of the Czech nation.  

Rychl\'\i{}k spent the rest of his life, especially the period 1953--58 when his 
retirement pension was drastically reduced, in poverty. The only activities he was allowed were 
connected with the work in the history of mathematics, especially 
treating the work of Bernard Bolzano. Apart from that Rychl\'\i{}k translated four books 
of Soviet mathematicians into Czech, the first one being Glivenko's textbook on probability theory. 
He died on 28 May 1968.  

Hyk\v sov\'a divided the scientific work of Karel Rychl\'\i{}k into five groups: algebra and
number theory (22 works), mathematical analysis (7), textbooks, popularization works, 
translations (16), works on Bernard Bolzano (14), and other works on history of 
mathematics (29). The most significant part is algebra and number theory in which Rychl\'\i{}k 
dealt with topical questions of building modern algebra. Most of his articles remained 
unnoticed in the world because they were published in national journals. 
The attention was drawn mainly to the paper 
{\it Zur Bewertungstheorie der algebraischen K\"orper} from 
1923 which was published in {\it Journal f\"ur die reine und angewandte
Mathematik} and was cited many times.

As far as mathematical analysis is concerned Rychl\'\i{}k published several works 
bordering algebra and number theory.  He was one of the first mathematicians to deal with  
$p$-adic analysis. 

Rychl\'\i{}k is an author of three textbooks the most important of which is 
{\it An Introduction to Elementary Number Theory} \cite{Rychlik1931}, 1931 which replaced 
the more than fifty-year old textbook of Studni\v cka. It was published in 1950 in a second 
revised and enlarged edition. Rychl\'\i{}k's textbook on probability theory \cite{Rychlik1938} is treated 
in Part 3 of this paper. The third Rychl\'\i{}k's textbook 
\cite{Rychlik1957} from 1957 deals with the methods of numerical solution of 
algebraic equations with real coefficients and it is based on his university lectures. 

From the point of view of Czech mathematics, Rychl\'\i{}k was famous, apart from his textbook 
on elementary number theory, for his interest in history of mathematics. 
Except for several popularizing articles about Czech as well as foreign 
mathematicians,\footnote{Let us recall at least his works concerning 
Abel's and Cauchy's stays in the Czech lands.} Rychl\'\i{}k made a name with his 
scientific work on part of Bolzano's inheritance. 
His interest in this area can be already seen in the 1920's when he published 
a paper dealing with Bolzano's continuous function which has a derivative at no point. 
Rychl\'\i{}k was a member of the Bolzano Committee under the Royal Bohemian Society of 
Sciences from 1924 and he took part in publishing Bolzano's works 
{\it Functionenlehre} (1930) and {\it Zahlentheorie} (1931). 
He continued to be a member of the Committee after the war, until 1952 
when it was suppressed. It was, however, restored in 1958 under the Czechoslovak 
Academy of Sciences, though it existed only briefly, until 1961. 
During this period no large Bolzano's works were published and Rychl\'\i{}k concentrated on 
writing comments on Bolzano's works connected to the theory of real numbers and logic. 

The work of Rychl\'\i{}k's on Bolzano's inheritance was famous internationally, 
he read a lecture {\it La th\'eorie des fonctions de
Bolzano} at the international mathematical congress in Bologna in 
1928.\footnote{Karel Rychl\'\i{}k took part only in two international mathematical 
congresses, apart from the one in Bologna, he was in Strasburg, 1920.} 
After the war Rychl\'\i{}k corresponded with some significant Bolzano researchers. 

Karel Rychl\'\i{}k was an extraordinary member of the Royal Bohemian Society of 
Sciences from 1922 and an extraordinary member of the Czech Academy of Sciences and Arts 
from 1924.

\section{Rychl\'\i{}k's textbook}

Rychl\'\i{}k's textbook \cite{Rychlik1938} is a booklet containing 144 pages. 
 It is not printed but was produced using 
a reprographical method by zincography, an ancestor of the offset printing where a `photography' of a text was obtained
using the static electric properties of paper covered with zinc oxydium.  
The booklet presents therefore  itself as a typed manuscript. All the special signs are hand-written, in particular 
all the mathematical signs and formulae: this tedious work had been realized by Otokar Pankraz, who is acknowledged 
for that in the foreword. Pankraz was Rychl\'\i{}k's assistant at the Technical University and we shall return to him 
at length in section \ref{section4}. The text is divided in chapters (with a main title), and in numbered sections and 
sub-sections from 1 to 47, immediately followed by complementary sections about set theory  numbered from 101 to 104. 
This strange gap in numbering may be an indication of the author's intention to complete and improve the text which is 
in the present state only a kind of first draft of a future book. With this respect, the text has the exact aspect of 
lecture notes written by a teacher for students studying probability theory. And the title itself where the word 
{\it \'Uvod} (Introduction) is highlighted contributes to the fact that the book is intended for beginners or at last 
for undergraduate students, who possess only the general mathematical bases, mostly in Analysis and Calculus, such as 
general properties of real functions and  basic (Riemann) integration.  The book contains a rather large bibliography, 
and also an index of the terminology, two properties not so common at the time (certainly not in textbooks devoted to 
students), which at least proves the author's remarkable modern pedagogical concern.

Although the book was published as a textbook for technical university students, and both Hyk\v sov\'a \cite{Hyksova2003} 
and Ma\v c\'ak \cite{Macak2005} state so, there might be arguments that indicate that the text was not intended primarily  
for technical universities. First, the exposition is not usual for technical university students, in the sense that it is 
too abstract, e.g. as far as we know the set theoretical conception of mathematics did not appear at technical 
universities in Czechoslovakia before WWII at all. Rychl\'\i{}k's option for this kind of presentation can be explained either by 
his highly innovative approach, or by the fact that he had a different audience in mind. The latter reason 
can be explained by another argument: Rychl\'\i{}k wanted to use the text for his lectures at Charles University, however, 
he used the opportunity to publish it with the support of the Central Publishing Committee of the Czech Technical University. 

The author states in the Foreword of the book that his aim is to compare axiomatic theory with other possible 
presentations of probability, nevertheless, he mainly concentrates on the axiomatic presentation. 
Let us therefore first compare Rychl\'\i{}k's title to the two other textbooks where Kolmogorov's axiomatic definition 
of probability was exposed before, namely Kolmogorov's own original text \cite{Kolmogorov1933} and 
Cramer's textbook \cite{Cramer1937}. In Kolmogorov's title, the world {\it Grundbegriffe} (foundations) is a clear 
allusion to the word  {\it Grundlagen} in Hilbert's famous treatise on the axiomatic foundations of geometry, and as 
such should have been understood by the mathematicians of the time as not devoted to students' use. This was also of 
course obvious from Kolmogorov's choice of Springer's collection {\it Ergebnisse der Mathematik und ihrer Grenzgebiete} 
for publishing his memoir, a collection that collected the most modern aspects of mathematics in progress and was clearly 
intended for specialists only.  As for Cramer's book, its title is clearly much more technical as it begins  by the words 
{\it random variables} and therefore also could certainly not reasonably be formulated for undergraduate students.

Therefore, as Hyk\v sov\'a \cite{Hyksova2003} had already mentioned, Rychl\'\i{}k's textbook seems to be a quite rare 
document about an attempt to teach an axiomatic version of probability theory to undergraduate technical university 
students in the 1930's. 
Rychl\'\i{}k chose a rather abrupt 
axiomatic approach of probability, using set theory, which could justify the argument that the textbook was intended for 
his Prague University students as technical university students were not familiar with even basic aspects of modern 
set theory. 
The author added a special part at 
the end of the booklet (from page 114) called {\it Mno\v ziny, Funkce, T\v elesa mno\v zinov\'e} which is to say Sets, 
functions, set fields (the name used for Boolean set algebras). In the first paragraph, he presents the basic properties 
of sets, beginning by a short presentation of the concept following Cantor's ideas. And indeed the section begins with 
a quotation of Cantor's assertions about what a set is from his paper \cite{Cantor1895}. The same quotation appears in 
the subsequent Pankraz's paper \cite{Pankraz1939} that we shall decribe in the next section.

He also refers to the new \v Cech's book on set theory \cite{Cech1936} advising the interested reader to consult 
it and the bibliography it contains. \v Cech was then professor of mathematics in Brno (where he had replaced Maty\'a\v s 
Lerch, a world famous specialist in number theory, after his death in 1922), and a prominent specialist in topology.  
As Frol\'\i{}k states in the Foreword to the third edition of Cech's book from 1974 (this edition is enlarged by the chapters 
which were found after \v Cech's death in the 1960's), this book played a very important role in introducing modern set theory 
in Czechoslovakia. Before there had not been even well-established Czech terminology concerning set theory. 
The book is amazingly dense, even for today's reader, and contains almost only the technical aspects 
of the subject, here again in a style which may certainly remind of Bourbaki's philosophy of mathematics. In its first 
edition, it contained four chapters, the last one devoted to abstract measure theory, where \v Cech introduced algebras 
and $\sigma$-algebras and measures built on them.  In particular, he builds Lebesgue measure over $\R^n$. This part on 
measure theory disappeared in the third edition of \v Cech's book. 

Rychl\'\i{}k insists on different operators on sets he would use continuously: inclusion, union and intersection 
(denoted as the sum $A+B$ and the product $AB$ of sets $A$ and $B$ as was usual). The paragraph about functions 
(pp.118-121) is quite amazing for a textbook written for students non specialized in mathematics. 
Rychl\'\i{}k gives the set-theoretical definition of function from  set $\UU $ to set  $\VV$ as a subset of the 
graph $\UU \times \VV$. Certainly, he would not have expected his students to read or at least to assimilate these 
notions. The fact that this part is situated in the last pages of the book may be a testimony of this. Nevertheless, 
he may have been intellectually satisfied to show the possibility of an entirely axiomatic construction based on set 
theory, a point of view selected in the same years by the founders of the Bourbaki group for the redaction of 
their {\it El\'ements de Math\'ematiques}. This also allows him to make an allusion on the notion of cardinality, 
defining countable and uncountable sets (pp.121-122), a difference which would play a role in his exposition. 
This section is followed by the definition of complement sets and a  presentation of de Morgan's algebraic manipulation 
of sets, a paragraph on combinatorics, and one about the {\it characteristic }(indicator) function of a set. 
Then comes the final section numbered 104 called {\it T\v elesa} (fields). Rychl\'\i{}k defines Boolean algebras 
on $E$, as a set $\AA$ of subsets of $E$ which is stable by sum and difference, and gives elementary examples. 
As he is quite exclusively concerned by probability theory in countable sets (and mostly in finite sets), he can 
reasonably limit himself to this situation. It is however not completely true as he also wants to give some 
considerations on real random variables. He therefore adds a last paragraph (104,4) where he defines a Boolean algebra 
on the interval $[A,B]$ as containing the empty set and all the sets which can be written as a disjoint sum 
$I_1+I_2+\dots +I_n+\{ p_1\} +\dots +\{ p_k\} $ where the $I_j$ are intervals (of any kind) included in $[A,B]$.  
He  gives a correct sketch of the proof that this is indeed a Boolean algebra, which he denotes by $\AA (A\dots B)$.  
It may be seem strange at first glance that Rychl\'\i{}k would not consider the usual situation (considered by 
Kolmogorov in his section 3) of the Boolean algebra $\BB_0$ of finite union of disjoint intervals of the type $[a,b]$. 
However, $\BB_0$ does not contain single sets $\{ p \}$ for which Rychl\'\i{}k would like to estimate probability when 
he considers real-valued random variables. Of course, Kolmogorov had no problem with this, as he immediately extends 
$\BB_0$ to the generated $\sigma$-algebra.

Now let us go back to the beginning of the book, which opens by a historical introduction. More precisely, this part 
is mostly concerned with philosophical considerations about the so-called classical definition of probability 
(the Laplacian definition as the number of favorable cases divided by the total number of cases) and the problems 
it generates, as well as the attempts by Cournot, Ellis and Venn to define only an {\it a posteriori} probability. 
Rychl\'\i{}k also mentions von Mises' construction of admissible sequences ({\it Kollektiv} - collectives: the word 
seems, however, not to be written down in the book) having correct limit relative frequences, and he refers to 
Copeland's paper \cite{Copeland1928}, as well as to Kamke's book \cite{Kamke1932}. Hyk\v sov\'a \cite{Hyksova2003} 
had already observed that Kamke's book seems to have had a great influence on Rychl\'\i{}k, which he had reviewed for 
the {\it \v Casopis pro p\v estov\'an\'\i{} matematiky a fysiky} and which had given him the occasion of his only research 
papers approaching probability theory. Moreover, after the appearance of Mises' book \cite{Mises1931} in 1931, 
Rychl\'\i{}k read a lecture on probability theory based on Mises' considerations at Charles University in Prague.

At the end of the introduction to his booklet, Rychl\'\i{}k writes that {\it there are various ways of determining 
probability. All these ways are attempted to be included also in the axiomatic method by means of laying down certain 
theorems which we do not prove - axioms - and other theorems are deduced as their logical consequence. Using this method 
is common in some parts of mathematics; the first attempts for axiomatization of probability theory were given by 
Bohlmann and Broggi.} The names of Bohlmann and Broggi belong to the long tradition which in Germany took over Hilbert's 
sixth problem of axiomatization of probability and appear in particular in the text presented by Bohlmann in 1908 at 
the International Congress in Rome. This story is deeply studied in \cite{ShaferVovk2005}. Finally, Rychl\'\i{}k 
concludes by asserting that he {\it will use Kolmogorov's axiomatization. Its usefulness is in expressing the results of 
the theory of probability by means of mathematical results: a random event is defined as a set of elementary events and 
the probability as a (real) set function.} One can observe that, basically, the justification given by Rychl\'\i{}k for 
the use of Kolmogorov's axioms is more or less the same as the one quoted by Cramer in \cite{Cramer1937}: it is a pure 
question of mathematical convention. Or, to say it more dramatically: {\it there is no justification at all}. As Cramer 
writes down: {\it the question of the [convergence of relative frequencies to the probability] will not at all enter into 
the mathematical theory}.  Though the whole text is written at non very high level, it seems that Rychl\'\i{}k wanted to 
convince his reader that probability theory is a part of mathematics. It is usual, when teaching at this level, to have 
rather general considerations about randomness and random experiences, and above all to use such considerations to 
generate intuition in the reader's mind (and this is indeed the specificity of probability calculus inside mathematics). 
On the contrary, Rychl\'\i{}k tries to get rid of such aspects any time it is possible as was of course the case with 
Kolmogorov's own book. However, as already mentioned, the Soviet mathematician had another kind of audience in view. A remarkable 
example of the fact is situated in the chapter devoted to independence. To give an example of three events $A,B,C$ such 
that $A,B$, $A,C$ and $B,C$ are independent, but $A,B,C$ are not, Rychl\'\i{}k repeats exactly Kolmogorov's redaction 
(page 10 of the {\it Grundbegriffe}, note 3) on a completely formal mode, without expressing a {\it concrete} situation 
-- e.g. with dice -- though it generally makes the result obvious to students. Poincare for example in \cite{Poincare1912} 
does not give other formulations than through card games, dice or urns models. Rychl\'\i{}k himself quotes Bohlmann's 
example from \cite{Bohlmann1908}, which is described with help of numbered black and white balls in an urn.

The proper lecture notes open by 5 pages introducing the vocabulary of probability. This is an occasion for Rychl\'\i{}k 
to furbish the (today!)  usual dictionary between set theory and probability folklore (set/event, 
disjoint sets/incompatible events and so on) which take place on pp. 8-9 and which he had directly copied from 
Kolmogorov's \cite{Kolmogorov1933} section 3 {\it Terminologische Vorbemerkungen}.

The next part is called {\it Klasick\'a definice pravd\v epodobnosti} (Classical definition of probability): 
Rychl\'\i{}k recalls Laplace definition and illustrates it by considerations on {\it Heads and Tails} and on dice. 
He immediately mentions  (page 13, paragraph 6) that the definition is 
unsatisfactory. 
To illustrate this, he recalls the classical d'Alembert's `mistake' where he obtained the value 2/3 for the probability 
of getting at least one {\it head} in two throws of a coin. D'Alembert claimed that it was sufficient to consider the set 
of events H,TH,TT as one may stop to play after having obtained one {\it head}. In fact, Laplace in \cite{Laplace1814} 
had already commented on d'Alembert's error and mentioned the necessary hypothesis that the elementary cases had to be 
supposed equally probable (and henceforth the definition of probability becomes circular\dots ). 
It can be also pointed out that placing an example of wrong probability calculations immediately after Laplace's classical 
definition of probability was usual -- see e.g. L\'aska \cite{Laska1921a}, pp.7--8 and Hostinsk\'y \cite{Hostinsky1950}, p.8, 
mention 
the example with throwing two dice, Kamke \cite{Kamke1932}, p.15 gives d'Alembert's error, and Czuber 
\cite{Czuber1900}, p.738 states both of them.  
Nevertheless, as  
\cite{Macak2005} suggests, Rychl\'\i{}k may have considered d'Alembert's mistake so convincing that he placed it 
before introducing his axiomatic presentation. 

Rychl\'\i{}k introduces his axiomatic presentation for a probability on a 
field $\AA$ with the axioms that are in the same order as in  the first section {\it Axiome} of Kolmogorov's booklet 
\cite{Kolmogorov1933}.  Like Kolmogorov, he gives as the first example the construction of probability on a finite set 
$E$ by attributing a nonnegative number $p_i$ to any element of the set $E$ such that the $p_i$ add to 1.

The following chapter, called {\it Zobrazen\'\i{} a ekvivalence pokus\accent23u a rozlo\v zen\'\i{} pravd\v epodobnosti} 
(Experiments mappings and equivalences and probability distributions) is an exposition of the transfer theorem. 
The presence of this rather theoretical chapter at this early stage seems mostly to offer a solution to the paradox 
mentioned above in the so-called classical definition of probability. Having mentioned the paradox, on page 13, 
Rychl\'\i{}k promised that in section 12,1\footnote{On page 13, the marked reference is section 21,1, which does not 
exist. There has clearly been an inversion in the figures.}, he would provide material to solve it. The aim of the chapter  
is to emphasize the fact that for two random experiments to be stochastically equivalent, it is not sufficient that the probability 
spaces containing the elementary events be equivalent which is to say that there is a one-to-one transformation between 
them. In the paradox case, the probability spaces are the same. One has also to look at the behaviour of the probabilities 
on the algebras on each of the sets. Therefore, one defines the stochastic equivalence in the following way: the 
probability spaces $(E,\AA ,P)$ and $(E',\AA ', P')$ are equivalent if there is an application $\varphi : E \rightarrow 
E'$ such that $P'(A')=P(A)$ if $A'=\varphi (A)$. Though Rychl\'\i{}k's presentation is a bit intricate, it is worth 
noticing that he tries to attract his students' attention to the basic modern idea of probability calculus: one has to get 
rid of the probability space and to work with distributions. It is remarkable that this idea comes so soon in his textbook, 
before exposing the classical results and tools of the theory or presenting any application. We however do not know how 
the students of the Technical university may have received the fact.

On the contrary, the two chapters to follow (pages 27 and 37, respectively) are ``expected'': the first one treats the 
(elementary) conditional probability, and the second one that of independence. The conditional probability of an event 
$A$ given an event $F$ such that $P(F)>0$ is defined as $\displaystyle {P(AF)\over P(F)}$, and, apart from the fact that 
Rychl\'\i{}k observes that it satisfies the axiomatic definition of probability, the chapter contains the usual facts: 
computations in the case of uniform probability, Bayes formula and the classical illustration with urns models.  

An interesting fact may be observed on page 33, where Rychl\'\i{}k presents Bayes' formula. He mentions the usual term 
{\it cause of $B$} for an event $A$ for which $P(B/A)$ is given, but also adds that {\it for Fr\'echet, this event should be 
preferably called the {\rm hypothesis} for the event $B$}. This sentence seems to be connected with Fr\'echet and 
Hallbwachs's book \cite{FrechetHalbwachs1924}, quoted in the bibliography. The title of  Chapter III of this book is 
{\it Probability of hypotheses (or causes)}. Fr\'echet, keeping a statistical point of view on mind, sees the core of 
Bayesian method as the fixation of a priori probabilities. For instance, in the subsection {\it Precautions to take when 
using Bayes' formula}, Fr\'echet insists on problems appearing {\it when the probability $\pi$ of different 
{\bf hypotheses } before the event are insufficiently known} (our emphasis). And he subsequently develops Poincar\'e's 
example of determination of the status of a card player, cheat or fair player, following the type of card he presents.

The long chapter on independence (pp.37-56) contains rather tedious considerations on the independence of sets. 
Rychl\'\i{}k then applies them to introduce the binomial distribution and to check its general properties.

The next chapter is devoted to the notion of expectation of a random variable. Though he could have defined a random 
variable with  more generality, Rychl\'\i{}k  chose for the moment to restrict himself to variables taking only a  
finite number of real values, so that the definition of expectation is immediate, but only through the distribution, as 
Rychl\'\i{}k is not able to define it as an integral. He therefore cannot avoid the usual ambiguity of the axiomatic 
presentation for elementary probability. For example, to prove that the expectation is linear (more precisely, that the 
expectation of the sum of two random variables is the sum of expectations) costs him the usual tedious proof - deported  
later in 33 pages 77-79 in the chapter about the law of large numbers where he cannot avoid it.  However, it is worth 
noticing that he makes use of the  notions of equivalence and transformation introduced before to explain that all the 
important facts about a (finite valued) random variable may be explained over a finite space (see in particular the 
paragraph at the top of page 59).

But it is the next chapter, p.63, called {\it Spojit\'e rozlo\v zen\'\i{} na p\v r\'\i{}mce a rozlo\v zen\'\i{} 
spo\v cetn\'e} (Continuous probability on the real line and countable probabilities) when Rychl\'\i{}k must try to do 
something with his axioms, though he lacks the structure of $\sigma$-fields. It must be mentioned that he mostly gets out 
of the trap with elegance. As Rychl\'\i{}k does want to present continuous probabilities also in an axiomatic way, he 
makes use of the already mentioned Boolean algebra $\AA (A\dots B)$ on $[A,B]$ composed by sets of the form
$I_1+I_2+\dots +I_n+\{ p_1\} +\dots +\{ p_k \}$. If one considers a non-negative, non-decreasing function $F$ such that 
$F(A)=0$ and $F(B)=1$ (Rychl\'\i{}k uses the same, most confusing, notation $P$ for the function $F$ and the probability 
built on it), a probability on $[A,B]$ equipped with the mentioned algebra is thus defined by the following properties
$$ P(\{ p\} )=0, P(I)=F(b)-F(a) \mbox{ if } I= (a,b).$$
Next (p.64), he asserts that {\it very often, the probability  on the field $\AA (A\dots B)$ is defined through
$$F(\alpha )=\int_A^\alpha p(t)dt$$ where $p$ is a non-negative, (Riemann) integrable function on every finite interval 
in $[A,B]$ and such that $\int_A^B p(u)du=1$}. As he is not able to properly define a real valued random variable, 
Rychl\'\i{}k defines only the {\it mathematical expectation of the function $f(\alpha )$} by the formula
$$E(f(\alpha ))=\int_A^Bf(\alpha )p(\alpha )d\alpha$$ {\it if the integral exists.} He then introduces Gaussian and 
uniform distributions through their densities.  The following paragraph (number 29) is devoted to the special case of a 
countable set $E=\{ \xi_1,\xi_2,\dots \}$, where the probability is as usual taken as $P( \xi_i)=p_i$ through a 
non-negative sequence $(p_i)$ adding to 1. He takes care of clarifying that the Boolean algebra in this case is composed 
by all the subsets of $E$. This allows him to define a (real) random variable $\alpha$ on $E$, simply as a function on 
$E$ taking values in $\R$. The mathematical expectation $E(\alpha )$ is then defined as the series 
$x_1p_1+\dots +x_np_n+\dots $ {\it if the series is absolutely convergent}. The justification given for this claim of 
absolute convergence is that it allows the result not to be ordering dependent.   In \cite{Kolmogorov1933} (p.33), 
Kolmogorov defines the expectation as an integral in {\it the sense of Fr\'echet}: for $x$ a random variable,  
if for any $\lambda >0$ the sum
$$S_\lambda =\sum_{-\infty}^\infty k\lambda P(k\lambda \le x < (k+1)\lambda )$$ is {\it absolutely} convergent and 
$S_\lambda$ converges to a limit when $\lambda \rightarrow 0$, this limit is the expectation of $x$. Therefore, 
integrability is of course equivalent to that of the modulus by construction. Kolmogorov does not mention the term order 
independence in the case of countable random values.

Maybe when he arrived at this point, Rychl\'\i{}k had  the impression that he must tell his reader something about a more 
general situation, and to make him feel that this would require more technicalities which are not of the current level of 
the booklet. He adds an {\it Observation} (paragraph 29,3 p.72) where he formulates the continuity axiom for the 
probability $P$ in the form of Kolmogorov's {\it Stetigkeitsaxiom} (p.13 of  \cite{Kolmogorov1933}): if $(A_n)$ is a 
non-increasing sequence in $\AA$ such that $\displaystyle \prod_{n=1} A_n=0$ then $\displaystyle \lim_{n\rightarrow +\infty}P(A_n)=0$. 
Rychl\'\i{}k does not mention the fact that in a Boolean algebra it may  occur that the set 
$\displaystyle \prod_{n=1} A_n=0$ were not in $\AA$. He proves that if the algebra $\AA$ is finite, the axiom is 
straightforward. And adds that the continuity axiom is also {\it valid for the countable case. In the case where the set 
of the elementary events is a Euclidean space the axiom [\dots ] follows from the properties of Lebesgue measure} for 
which he refers to the construction of Lebesgue measure in the aforementioned fourth chapter of \v Cech's book on set 
theory \cite{Cech1936}. It is interesting to observe however that \v Cech treats the case of general abstract measures 
in the book,  which satisfy the set continuity property (19.2.3 p.137) and not only the case of Lebesgue measure. 
Maybe Rychl\'\i{}k found too abstract for his students to imagine another measure than Lebesgue measure.

 The last chapter is called {\it Posloupnostn\'\i{} model pro rozlo\v zen\'\i{} pravd\v epodobnosti} 
(Sequential model of probability spaces). It is an attempt to make the junction between the axiomatic point of view and 
sequential definition of probability in the style of von Mises's model of {\it Collectives}. The reader of Rychl\'\i{}k's 
text may feel uncomfortable with the rather abrupt transition with the rather linear exposition of axiomatic method 
followed by the author since the beginning, however, Rychl\'\i{}k points out in the preface that the book concerns 
also the relationship to older mathematical theories of probability. 
It is reasonable to think that Rychl\'\i{}k may have felt necessary to try to justify that the axiomatic definition did 
not prevent the common use of probability theory based on frequences (such as in statistical classical models, such as 
mortality tables he gives as an example in section 46). As Shafer and Vovk had already observed \cite{ShaferVovk2005}, 
the German mathematical tradition (to which Czech scholars were still so close) made a clear distinction between 
theoretical probability (dealt with by mathematicians) and questions about its application discussed by experimental 
scientists and philosophers.  Urban's book \cite{Urban1923} (mentioned in Rychl\'\i{}k's bibliography) for example gives 
an interesting example presenting  both aspects, but one after the other.
 Urban was an amazing universal mind, born in Br\"unn-Brno in Moravia in 1884 and whose thrilling life had espoused 
the shaken history of his country: his biography is narrated in \cite{Ertle1977}. See also Bru \cite{Bru2003}. In 1923, 
he published the aforementioned book where he sums up his views on probability.  The first part of the book 
(chapters I to III) deals with randomness and probability on a philosophical point. Urban proves there to have an 
extremely good knowledge of the literature on probability and chance and in three cultural traditions 
(English, French and German) which seems to be an exception at the time. Certainly, his situation of member of a 
German community, in the Czech land  and having been for 7 years in Philadelphia (USA) may at least partly explain this 
very large overview.   The slightly smaller second part of the book (chapters IV and V) deals with the mathematical 
theory of probability. At first glance, it seems that there is no  connection between the two sides of the book apart 
from the fact that the mathematical probability of an event is a real number attributed to this event, with rules inspired 
by our views on randomness. Urban writes in the introduction of chapter IV: 

\begin{quote}
{\it The statements of the calculus of 
probability are abstract and are not propositions about real events. Their logic is internal, and nothing tells whether 
there exist objects corresponding to the conditions declared in these propositions. Neither with Bernoulli's statements 
nor with Poisson's theorem or any other more deep examination,  it is possible to describe the reality. Every use of the 
calculus of probability in concrete situations must be preceded by an examination to check whether the  studied processes 
and phenomena have the properties required by the theory.}
\end{quote}

 As alreasy mentioned, Kolmogorov in the {\it Grundbegriffe} had evacuated the problem since the very beginning: he asserts that the question 
of the concrete interpretation of probability is not a mathematical question. To explain why Kolmogorov is so brief on 
the subject, it should also be noticed that he had already discussed this point before, and this time at length, in the 
introduction of his perhaps most important work in probability theory where he introduced continuous Markov processes 
\cite{Kolmogorov1931}.  However, Rychl\'\i{}k's purpose was to teach students,  non specializing in mathematics, who 
moreover were supposed to be in contact with concrete financial application such as actuarial mathematics. Therefore, it was 
certainly necessary to say something.
 It is possible that the part about the sequential model is directly inherited from his former teaching of probability 
calculus following von Mises axiomatization that he prepared for Charles University in academic year 1931-32. In sections 
37 to 41, he presents models of relative frequency of an event in a sequence, and shows that under reasonable hypotheses, 
the relative frequency limit satisfies the general properties required for probability. Rychl\'\i{}k develops the 
classical argument for foundation of probability theory on properties of relative frequencies. In Kolmogorov's booklet,  
the argument is briefly mentioned as {\it Empirische Deduktion der Axiome} p.4. The interesting part is section 42, 
called {\it Posloupnostn\'\i{} model pro rozlo\v zen\'\i{} pravd\v epodobnosti} (Sequential model of probability 
distribution). There he intends to join the axiomatic definition and the relative frequency model by defining a Boolean 
algebra as the collection of subsets for which the limit of relative frequencies exists. His main reference is Kamke's 
book (\cite{Kamke1932})  where Kamke formulates a (rather intricate) sequential model to define what he calls a 
{\it Wahrscheinlichkeitsfeld} as a sequence of sequences having good  relative frequency properties. As said earlier, 
Hyk\v sov\'a \cite{Hyksova2003} mentioned the importance of Kamke's book for Rychl\'\i{}k who seems to have been very 
convinced by the model. He wrote an enthusiastic review of the book in the {\it \v Casopis} (\cite{Rychlik1932}). 
He mentioned there:

\begin{quote}  
{\it  while at the origin of the theory of relativity, the geometry of Euclidean and 
non-Euclidean spaces of three and more dimensions, which physicists needed, were already developed on a high level, 
mathematics are not in such a perfect state as far as probability theory is concerned. There are serious objections 
against the way how probability theory is presented nowadays. There is still much to be done to have a perfect building 
of the theory and the consequences of the induced theorems are very often overestimated. However, as probability theory 
has gained a great importance, it is necessary to build it in the same clear and precise way as geometry is built. 
Nobody denies this need as far as presenting geometry is concerned. }
\end{quote}

Rychl\'\i{}k produced soon afterwards himself a paper on well chosen binary sequences \cite{Rychlik1933}: in the presently 
discussed chapter of Rychl\'\i{}k's booklet, the section 39 is apart from others devoted to the case of binary sequences 
of 0 and 1. The chapter ends by considerations about possible application of the sequential probability model, 
in particular with the already mentioned section (number 46) about mortality tables which may be an appealing example 
for students whose destiny should be to work in insurance companies.

\section{Reception and destiny of Rychl\'\i{}k's booklet}\label{section4}

It would certainly be extremely interesting to recover some impressions the students may have had by listening to 
Rychl\'\i{}k's lectures, as well as to know more precisely what he actually taught during these lectures. Unfortunately, 
it seems that we lack data for these facts.  
The main trace which is left is constituted by Pankraz's comments. Otomar Pankraz (1903-1976), already mentioned, 
was Rychl\'\i{}k's assistant at the Technical University.

Otomar Pankraz was born on 25 March 1903 in Nov\'e Dvory near P\'\i{}sek.\footnote{The information 
about Pankraz's life were gathered from his personal files kept in the archives of 
Charles University, Prague Technical University and in the fund of the Ministry of 
Education in the National Archives in Prague.} His secondary school education was negatively
influenced by a poor material situation of the family, which resulted in the fact that Pankraz
could not finish a comprehensive secondary school, but graduated from a secondary technical school in Prague in 1923. 
This secondary school graduation exam did not allow him to attend a university as an ordinary student, 
he, therefore, followed the lectures at Charles University only as an extraordinary student in years 1923--1929. Only 
a supplementary secondary school graduation examination from 1929 enabled him to apply for doctoral examinations 
from mathematical analysis, algebra and theoretical physics. He obtained his doctoral degree 
in natural sciences in 1931. Although he had also fulfilled all the requirements needed for a future 
secondary school mathematics and physics teacher, he never took the final teacher examination. 

In June 1929 Pankraz started working as an actuarial mathematician for preparing superannuation 
scheme in the General Pension Institute. However, it was only in 1931 that he passed actuarial 
mathematics examination and mathematical statistics examination with professor Emil Schoenbaum 
who worked as director of this institute in 1919--1939. Then Pankraz decided to follow an academic career 
with the possibility of scientific work and in May 1931 he became a mathematics assistant 
of Karel Rychl\'\i{}k at the Prague Czech University. He remained at this post until the closure of Czech 
universities and he was deprived of it officially only in March 1945. 

In 1935 Pankraz habilitated at the university with the thesis {\it Zur
Grundgleichung f\"ur den zeitlichen
Zerfall der statistischen Kollektivs} \cite{Pankraz1933c,Pankraz1933a} in actuarial 
mathematics and mathematical statistics. According to Schoenbaum the work {\it shows the author's complete 
mastering of analytical tools for solving complicated integral--differential equations.} Pankraz's 
scientific activities were also favourably commented on by Karel
Engli\v s, the governor of the National Czechoslovak Bank, and actuarial mathematics professors 
L\"owy (Heidelberg), Moser (Bern) and Cantelli (Rome). During his habilitation colloquium Pankraz was asked, apart from others, 
in great detail about Mises theory of probability. 

From the summer semester 1936 Pankraz gave regular lectures as a private docent at the university, with the 
exception of the summer semester 1937 when he went to London for a research fellowship financed by {\it Rockefeller
Foundations}. In the academic year 1938/39 he read lectures on probability theory instead of professor 
K\"ossler.

Based on the habilitation which took place in 1937 Pankraz was appointed a docent also at the 
Prague Technical University in 1938. 
He submitted two works on integral equations \cite{Pankraz1933b,Pankraz1936} as his habilitation work. 
Positive opinions were given by professors Rychl\'\i{}k a Hru\v ska who also praised other 
Pankraz's extensive publication activities. In 1938/39 Pankraz gave lectures on integral equations and their applications in 
engineering at the Technical University. 

After the closure of Czech universities in November 1939 Pankraz lived in Prague and he received his assistant's salary 
until January 1945 although he did not perform any activity.\footnote{During the war he translated the book by 
E.~Wagemann {\it Wo kommt das viele Geld her?} into Czech in 1943.} The payment of the salary was stopped after he had refused 
doing the work he was charged with by the Ministry of Education for more than a year. 

Pankraz was arrested after the war, already on 19 May, and on 21 June 1946 he was found guilty of 
{\it propagating and supporting the Nazi movement especially by approving murdering university teachers after the execution 
of Heydrich and closing the universities and by praising the Nazism and the work of Czech traitors for Germany} by 
an extraordinary people's court. The court also managed to prove that the attorney-in-law Dr. Chlubna was arrested 
on the basis of his information. The court sentenced Pakraz to five years of heavy jail, the loss of property and 10 years 
of the loss of citizen's honour which he served in special working groups. 

We know nothing about the later life of Otomar Pankraz except that he died in Prague on 12 December 1976 at the age of 73. 

Pankraz, like Rychl\'\i{}k, seems to have been deeply interested in the evolution of probability theory during the 1930's. 
He regularly made reviews in the {\it \v Casopis}. In particular, it was him who reviewed von Mises' book: 
in his review \cite{Pankraz1931}, 1931, he not only describes the content of the book but also (briefly) presents the 
objections opposed to von Mises about his {\t collectives} theory and takes Mises' side. He writes: 

\begin{quote}
{\it There are various objections against this theory. For example,  it is said that it is not allowed to use the 
analytical concept of limit  in probability theory. This objection is not valid because it is obvious from how Mises 
deduces the laws of large numbers that an analytical limit is absolutely satisfactory.  However, the objections based on 
the principle of impossibility of gambling system and on the combination operation are more serious. I can see a solution 
to these difficulties in stating that Mises' requirements are principles which are closely connected to experience, not 
purely logic axioms (though derived from experience)}.\footnote{For a discussion of early objections against Mises' theory, 
one may consult \cite{Vonplato1994}, pp.192-197.}
\end{quote}

Pankraz had to face Hostinsk\'y's unsatisfaction about one of his reviews. Pankraz wrote in \cite{Pankraz1932}, 1932,  the 
review of Hostinsk\'y's little treatise on Markov chains (his first really internationally celebrated publication) 
\cite{Hostinsky1931} 
and seems to have been deluded by the title {\it M\'ethodes g\'en\'erales} as he regrets that 
Hostinsk\'y had not  given a general exposition of probability theory. Maybe the ambiguity of Hostinsk\'y's title 
explains also why it was the (young) Pankraz who had been given the task of reviewing the text of an internationally 
known scientist as Hostinsk\'y. Despite the review is quite laudative, Hostinsk\'y seems to have been taken aback with it 
and subsequently wrote a rather dry answer, also published in the {\it \v Casopis} (\cite{Hostinsky1932}).

In 1938, Pankraz took over the task of reviewing Rychl\'\i{}k's booklet.  Though reviews of that time were often 
longer than today as the reviewer included his own opinions, the length of Pankraz's review  
\cite{Pankraz1938} is quite noticeable (two full pages and half: as a comparison the review of von Mises' book is only one 
page long), which at least proves a particular interest in the question (of course, self promotion might have been also 
involved!).  Pankraz justifies the necessity of a book by the quick developments of the theory of probability, and mentions 
that Rychl\'\i{}k had chosen to exploit Kolmogorov's axiomatization. He adds that {\it contrary to Kolmogorov who 
(\dots ) only briefly exposed the problems, professor Rychl\'\i{}k presents the axioms in detail}. This slightly 
exaggerated sentence can be understood as an indication that the levels of the two texts are very different. Then 
follows a description of Rychl\'\i{}k's textbook. Pankraz comments at length on the bases of the axiomatic foundation: 
\begin{quote}
{\it The readers who are used to think with the concepts of classical probability theory must be warned that the 
transition to the new theory lies exactly in the definition of the random event. In the classical theory, this was a 
vague notion and its justification was based more on intuition than on mathematical bases. The core of the new theory is 
precisely to state the most important feature that it is necessary to consider a random event as a collection (i.e. a set) 
of events, not as an isolated event. There is also a conceptual difference between an element and a collection (set) of 
elements which is met nowadays in nearly every branch of pure and applied mathematics (and also logic) which makes the 
difference between ``modern''  and ``classical'' theories. Naturally, this difference has been known for a long time but 
is properly considered only in modern theories. } 
\end{quote}
After presenting the axiomatic definition of probability on a Boolean 
algebra, Pankraz adds that in his book Rychl\'\i{}k {\it shows very clearly and in details that it is indeed possible to 
build probability theory on the basis of these axioms and that classical probability theory is included in it. The 
transition to applications, mainly in statistics, is enabled by the so-called sequential model of probability 
distribution.} Once again, Pankraz describes the sequential model in details. The end of the review expresses a real 
enthusiasm: 
\begin{quote}
{\it The set-theoretical considerations thus (1) provide us with an exact and simplified basis for 
building probability theory and (2) lead to an immediate application in statistics. At the same time, the results are 
more extensive than in the classical theory. The relevance of the new methods is indisputable and it now depends on 
didactical methods how it may be possible to transfer the way of thinking from the old form to the new one in the 
simplest way. It is only a question of schools and it is a good thing that the Czech Technical University in Prague is 
ready to support such new methods.

Rychl\'\i{}k's book, although it is modestly called an Introduction, presents author's original thoughts at numerous 
places and indicates directions in which the study of probability theory can be extended.} 
\end{quote}

One passage is rather remarkable 
for a review. Pankraz writes that {\it [T]he book is not printed, but realized with zincographic reproduction which 
enables the author to publish it later in a modified state which could reflect the new and latest results in probability 
theory. This method is much more suitable than publishing expensive printed books for the Czech environment where the 
sale of mathematical literature is rather small}. What Pankraz intends with this sentence is not quite clear. 
Apart from the reasonable mention of the economic side of the question, the sentence may also suggests that the two men 
had discussed future possible improvements of the booklet. At least, the sentence seems to assert that the present state 
of the text is not completely satisfactory and that the work is in progress. Another hint of this is the fact that in 
winter term 1938-39 and summer term 1939, Pankraz had taught lectures on probability theory at Charles University and 
had written his paper \cite{Pankraz1939} published in the Prague {\it Aktu\'arsk\'e v\v edy} (Journal of Actuars). 
The paper may  reflect the content of Pankraz's lectures, and is mostly devoted to the question of axiomatization of 
probability. For Pankraz, the two concurrent ways through which probability has been defined, i.e. Mises frequentist 
approach and Kolmogorov's axiomatic method generate unsolved problems. In Section 3 of his paper, he sums up these 
problems :  1) the  internal contradiction of Mises theory which pretends to give a precise definition of probability 
though it uses concepts of subjective nature. 2) the incompleteness of Kolmogorov's axiomatization to which Kolmogorov 
is obliged to add a separate definition for conditional probability.

 It is nevertheless clear in  the paper that the accent is mainly put on Kolmogorov's axiomatic construction 
(quite remarkably,  von Mises, though quoted in the text, is absent from the bibliography). Pankraz describes at length 
large bases of set theory in his section 5 (where, as already mentioned, he quotes Cantor's assertions on sets as 
Rychl\'\i{}k had done) and in section 6 connects these set-theoretical notions to the notion of random event. We had 
already observed in his review of Rychl\'\i{}k's book how he emphasized the fact that the set conception of the random 
event has become the basic notion of modern probability theory. The section 7 is entirely devoted to a brief description 
of Kolmogorov's axiomatics and to criticizing it in the light of the his previous comparison of Mises and Kolmogorov 
definitions.  For Pankraz, the introduction of conditional probability as a separate notion by Kolmogorov cannot be 
justified on a logical basis: 
\begin{quote}
{\it In fact, the question is whether probability is defined as a function of one or two 
arguments. Kolmogorov tries to get along with the set function $P(A)$ with one argument, defined through axioms. However, 
at the same time, he comes to the fact that this one-argument function is not enough for formulating Bayes formula, and 
therefore he introduces the set function  $$P_A(B)={P(AB)\over P(A)}$$ with two arguments $A$ and $B$, which he proves 
also to satisfy his axioms. He must, however, suppose that the set $A$ is fixed, and therefore this assumptions allows 
him to see the two-argument function as a one-argument function which obviously satisfies the axioms.

Reichenbach has already pointed out that this method is not satisfactory. While Reichenbach's objection to the 
introduction of the new symbol $P_A(B)$ was more formal than mathematical, we find the condition of $A$ being constant 
in the definition of $P_A(B)$ inadmissible from a mathematical point of view as in this case it would not be possible for 
$A$ and $B$ to be variable at the same time. }(\cite{Pankraz1939}, [7,2]). 
\end{quote}
Reading Pankraz's bibliography, it seems to 
have been quite interested by the discussions about the status of probability in quantum mechanics, which flourished 
during the first third of 20th century. They have been the study of numerous papers: see for example \cite{Vonplato1994}, 
Chapter 4. A famous international conference organized by the philosophers from Berlin and Vienna Circles about the 
implications of the new physical theories was held in Prague (as a barycenter between the two towns!)  in 1929. Here were 
present Carnap, Reichenbach, von Mises and many other of the names mentioned in the present paper.  As Hyk\v sov\'a has 
already briefly suggested in \cite{Hyksova2006}, it is quite possible that Rychl\'\i{}k and Pankraz have also attended 
it. Pankraz had certainly at least read the texts of the conferences published in the founding issue of the journal 
{\it  Erkenntnis}. In particular, he quotes Waismann's conference on the logical foundations of probability 
\cite{Waismann1930}.  The two last sections of \cite{Pankraz1939} are besides devoted to considerations about quantum 
mechanics.

The previous remark on Kolmogorov's axiomatization leads Pankraz to formulate his own system of axioms for a probability 
defined as a two-arguments set function, which he claims to be complete and non-contradictory. The probability is 
therefore defined as satisfying an additivity property on the second argument
$$P(A,B+C)=P(A,B)+P(A,C)$$ if $B$ and $C$ are disjoint sets (he mentions a version with $\sigma$-additivity on the next page), and the multiplication axiom as
$$P(A,BC)=P(A,B)P(AB,C)$$ when $P(A,B)>0$.  
It can be seen that Pankraz formulates the same kind of axiomatization as Popper would find necessary to  formulate only 
some twenty years later in the Appendix IV of \cite{Popper1959}.  
It should be observed that in  1938, Popper \cite{Popper1938} had proposed an alternative set of axioms where the 
primitive notion is the `absolute'  probability.  In the Appendix II of \cite{Popper1959}, where the beginning of  
\cite{Popper1938} is printed again, Popper asserts (with a little exaggeration) that he was the first to propose that 
the  mathematical theory of probability should be elaborated as a {\it formal system}. Popper's paper is not mentioned 
in Pankraz's text though it was published in the journal {\it Mind} which Pankraz happened to know as it is quoted for 
Nagel's extensive review on  Reichenbach's book in the bibliography of \cite{Pankraz1939}. Maybe, the 1938 issue of the 
journal had not yet reached him at the time. However, he insists that {\it to speak about the probability of the event 
$B$ without having mentioned the [reference] event $A$ has no meaning}. Pankraz therefore joins the trend of those who 
consider that the basic notion is conditional probability,  as estimating the probability of an event makes sense only 
with reference to another event: there is not such a thing as absolute probability. This leads some of the main 
representatives  of this trend to refuse an axiomatic foundation, as de Finetti and his entirely subjective conception 
of probability. Much later, de Finetti exposed the following opinion in the introduction to Chapter 4 of his textbook on 
probability theory: {\it Every prevision, and, in particular, every evaluation of probability, is conditional; not only 
on the mentality or psychology of the individual involved, at the time in question, but also, and especially, on the state 
of information in which he finds himself at that moment} \cite{Finetti1970}.  It should be extremely interesting to know 
what Rychl\'\i{}k or Pankraz, two seemingly  open minded scientists  who enthusiastically successively adopted von Mises 
and Kolmogorov's points of view on probability, may have thought of de Finetti's (being also an actuar like Pankraz was) 
views.
Up to now, we have no trace to which extend Pankraz knew de Finetti's texts (such as his comments 
on Bayes theorem in \cite{Finetti1939}, pp.20--23), the only place where he mentions 
his name is the review \cite{Pankraz1940c} of Nagel's book \cite{Nagel1939}.

The next year 1940, Pankraz published a new paper, this time in the {\it \v Casopis} \cite{Pankraz1940}, where he came 
back to the fundamental notions exposed in his 1939 paper, to satisfy  the claim of his students looking for a smooth 
introduction to the new aspects of probability theory, as he himself justifies in the foreword (\cite{Pankraz1940}, p.D73). 
Quite interestingly, the name of Glivenko appears in the paper, through the book \cite{Glivenko1938} published in Paris 
in 1938 by the Soviet mathematician about general algebraic structures where he mentioned some connection with the 
foundation of probability theory. In \cite{Glivenko1938}, p.26, Glivenko proposes to define the stochastic type of 
an event by considering that $A\sim B$ if and only if $P(A/B^c)=0$ and $P(B/A^c)=0$. On the {\it structure $S$ whose 
elements are events, the probability defines a norm} (p.26).  From Glivenko's book appears  for instance for the first 
time in Pankraz's paper the expression Boolean algebra (p. D79).  At that time, Glivenko had also written a textbook on 
probability theory \cite{Glivenko1939}, published in 1939 where he made use of Kolmogorov axiomatization. The book reached 
Rychl\'\i{}k somehow during WWII  and was translated by him into Czech; it was published in Czechoslovakia only in 
1950 (\cite{Glivenko1950}).
Basically, the paper presents also the construction of probability as a two-argument function. The text is however more 
strictly technical and Pankraz has left aside many of the philosophical questions he exposed in \cite{Pankraz1939}.  
His aims seems to show how it is possible to rapidly construct a probability function with his axioms, which is 
practicle to manipulate with.  Nevertheless, up to this point Pankraz does not seem to have realized the problem in his 
system, where $P(\emptyset , \emptyset )=1$ and where $X\subset Y$ may lead to probabilities greater than 1. In 
a subsequent note \cite{Pankraz1940b}, Pankraz proposed a corrected version of his axioms. Not without some irony, 
he mentions that  it was his reading of Reichenbach's book that lead him to this error. Naturally, one may reasonably 
comment that the obtained system of axioms, where it is now necessary to check the non-emptiness  of the considered sets 
for the relations on probability to be valid is not as comfortable as was his former (incorrect) one.  However, the 
paper \cite{Pankraz1940b} contains more than a simple technical correction, as it gives a separate status to the 
multiplication axiom. Pankraz's axioms are now the following ones

 (i) $\AA$  is an algebra of subsets of $E$.\footnote{Pankraz's original notation for the algebra is $\Omega$! 
$\AA$ seems preferable for our comments, in order to avoid the nowadays compulsory understanding of $\Omega$ as the 
set of elementary events.}

(ii) For every couple $(X,Y)$ of elements of $\AA$ such that $X\neq \emptyset$, $P(X,Y)$ is well defined and non-negative.

(iii) For every couple $(X,Y)$ of non-empty elements of $\AA$ such that $Y\supset X$, $P(X,Y)=1$.

(iv) For any $X,Y,Z$ in $\AA$ such that $X\neq \emptyset $ and $YZ=\emptyset$, one has $P(X,Y+Z)=P(X,Y)+P(X,Z).$

One may in particular observe that $P(X,Y)$, interpreted as {\it the probability of event $Y$ given that 
event $X$ occurs} is well defined even if $P(X)=0$ (unless $X=\emptyset$ of course). This problem of the conditioning 
on 0 probability set had been for a long time a recurrent question. A celebrated example is the {\it great circle paradox} 
(sometimes called the {\it Borel-Kolmogorov paradox} ) proposed by Bertrand among his famous list of paradoxes, and 
investigated by Borel in \cite{Borel1909}: if a point is chosen {\it at random} on a sphere, the point is {\it uniformly} 
chosen on the equator but {\it not uniformly} on any meridian. The investigations around Bertrand's paradox have been 
studied in \cite{ShaferVovk2005}.  Borel has pointed that the status of the conditioning proposition, and that the 
source of the ``paradox'' was in the fact that at first glance, the conditions {\it the point is on the equator} or 
{\it the point is on a meridian} have probability 0.  Kolmogorov, in \cite{Kolmogorov1933} (p. 44), mentions Borel's 
paradox and claims that his own definition of conditional distribution via Radon-Nikodym theorem gave a satisfactory 
explanation for the paradox. Nevertheless, the impossibility of conditioning by a null probability event had never been 
considered as completely satisfactory, a classical observation being that it is hard not to be able to assert that 
$P(X/X)=1$. Axiom (ii) of Pankraz gives indeed that $P(X,X)=1$ for any non empty set from $\AA$.

Pankraz observes that his axioms (i) to (iv) are complete and non contradictory. However, up to this point, nothing 
allows to make a connection with a statistical interpretation. As it was the case in Rychl\'\i{}k's textbook as we had 
already commented on, this possibility of embedding the classical frequency interpretation of probability into any 
proposed axiomatic definition seems to have been a major concern of the two Czech mathematicians. Pankraz comments: 
{\it It is necessary to agree on the meaning of the expression `statistical' interpretation. Generally one agrees that 
every number which can be reduced to frequencies is `statistical'. Indeed, apart from a single case (quantum physics) 
this agreement never leads to doubts.} He therefore proposes his fifth axiom:

(v) For every couple $(X,Y)$ of elements in $\AA$ such that $X\neq \emptyset$, $P(X,Y)=\displaystyle {P(E,XY)\over P(E,X)}$.

Pankraz should certainly preferably have written the axiom as $P(X,Y).P(E,X)=P(E,XY)$, to avoid a subsequent intricate 
justification of what happens if $P(E,X)=0$. Then, asserts Pankraz, $P(E,XY)=0$ by axiom (iv) and so the fraction is not 
defined and this seems to mean that axiom (v) is not concerned by this case. Anyway, in the case where $P(E,X)>0$, 
the quantity $P(X,Y)$ is indeed equal to the quotient $\displaystyle  {P(E,XY)\over P(E,X)}$. Pankraz does not in fact 
explain why this quotient may be interpreted as a frequency: he had probably in mind a sequential model of the same kind 
as Rychl\'\i{}k's that we have presented in the previous section. From the comments he adds, it seems that his most 
concern was about the possibility or impossibility of building $P(X,Y)$ from a common reference set, a question directly stemming 
from quantum physics:
\begin{quote}
 {\it The question of validity of axiom (v) for quantum physics remains undecided. Is it possible to interpret every case 
of this physics by means of frequencies if we consider that there exist

1. the least measurable length of order $10^{-13}$ cm

2. the least measurable time interval of digit place $10^{-13}$ cm/$c$ , $c$ = velocity of light in the vacuum, and

3. complementary measurable data following Heisenberg uncertainty principle?

Is it justified to add axiom V to axioms I -- IV for the `probabilistic' description of  subatomic events taking place 
in spatio-temporal fields of order smaller than $10^{-13}$ where, in general, I cannot speak about the possibility of 
`counting' elements though these are collective events? As it is obvious, the question is whether the probability 
exposition of specific events necessarily means the exposition by means of frequencies, or whether its framework is 
wider.}
\end{quote}

Pankraz concludes his paper by formulating the Bayes formula with his two-argument probability axiomatization.

\section*{Conclusion}

In March 1939 German troops entered Prague, and after few months of occupation, the Germans decided to close the 
universities. This moment marked also the end of the original  teaching experience by Rychl\'\i{}k and Pankraz. 
In the first place, as was narrated before, the two had stopped their mathematical activity soon after the war. In the 
second place, Kolmogorov's axiomatization gradually imposed itself as the most practicle presentation of probability 
theory and therefore became more and more adopted. In the 1960's, almost every mathematician involved in probability 
adopted it. Exceptions may nevertheless be quoted. The `old' probabilists, such as Paul L\'evy,  refused to change 
the presentation they were used to -a thing which seemingly did not prevent L\'evy  from brilliant discoveries on Brownian 
motion during the 1950's (see \cite{BarbutLockerMazliak2004})! In Italy also, Kolmogorov's axiomatics was slow to appear, 
due to the obstinate opposition of de Finetti and his great personal influence in the Italian community. And of course, 
we may mention the strong opposition in France to the abstract measure theory by the tenors of the Bourbaki group.  
However, in Czechoslovakia after the installation of a Soviet Union inspired governement, the situation was clear and 
Kolmogorov's axiomatization imposed itself as in the USSR. Rychl\'\i{}k even obtained some consideration from the new 
oriented politicians for having translated Glivenko's book from Russian before 1945.  As already 
mentioned, the publication of his textbook on axiomatic probability theory in 1938 had been a surprising event without real consequences.  But  we 
do not know what they may have been in a quieter context\dots

\end{document}